\documentclass[11pt]{article}
\usepackage[cm]{fullpage}
\usepackage{amscd}
\usepackage{amsmath}
\usepackage{amsfonts}

\newcommand{\ep}{\epsilon}

\newcommand{\dbar}{\overline{\partial}}
\newcommand{\ddbar}{\sqrt{-1} \partial \overline{\partial}}

\newcommand{\Rc}{\textrm{Ric}}

\newcommand{\op}{\omega_{\phi}}
\newcommand{\opsi}{\omega_{\psi}}
\newcommand{\ddbp}{\sqrt{-1}\partial \overline{\partial} \phi}
\newcommand{\ddbpsi}{\sqrt{-1}\partial \overline{\partial} \psi}
\newcommand{\dotp}{\dot{\phi}}
\newcommand{\dotpsi}{\dot{\psi}}
\newcommand{\AK}{\mathcal{K}^{\mbox{\tiny{-1}}}}
\newcommand{\kahler}{\mathbb{K}}
\newcommand{\kahlerp}{\mathbb{K}^+}

\begin{document}
\newcounter{remark}
\newcounter{theor}
\setcounter{remark}{0}
\setcounter{theor}{1}
\newtheorem{claim}{Claim}
\newtheorem{theorem}{Theorem}[section]
\newtheorem{proposition}{Proposition}[section]
\newtheorem{lemma}{Lemma}[section]
\newtheorem{defn}{Definition}[theor]
\newtheorem{corollary}{Corollary}[section]
\newenvironment{proof}[1][Proof]{\begin{trivlist}
\item[\hskip \labelsep {\bfseries #1}]}{\end{trivlist}}
\newenvironment{remark}[1][Remark]{\addtocounter{remark}{1} \begin{trivlist}
\item[\hskip
\labelsep {\bfseries #1  \thesection.\theremark}]}{\end{trivlist}}

\centerline{\bf ENERGY FUNCTIONALS AND CANONICAL K\"AHLER METRICS\footnote{The second author is supported in part by National Science Foundation grant DMS-05-04285} }

\bigskip
\begin{center}
\begin{tabular}{ccc}
{\bf Jian Song}  & & {\bf Ben Weinkove} \\
Johns Hopkins University & & Harvard University \\
Department of Mathematics & & Department of Mathematics \\
Baltimore MD 21218 & & Cambridge MA 02138 \\
\end{tabular}
\end{center}
\bigskip

\noindent
{\bf Abstract.} \ Yau conjectured that a Fano manifold admits a K\"ahler-Einstein metric if and only if it is stable in the sense of geometric invariant theory. There has been much progress on this conjecture by Tian, Donaldson and others. The Mabuchi energy functional plays a central role in these ideas.
We study the $E_k$ functionals introduced by X.X. Chen and G. Tian which generalize the Mabuchi energy. We show that if a Fano manifold admits a K\"ahler-Einstein metric then the functional $E_1$ is bounded from below, and, modulo holomorphic vector fields, is proper. This answers affirmatively a question raised by Chen. We show in fact that $E_1$ is proper if and only if there exists a K\"ahler-Einstein metric, giving a new analytic criterion for the existence of this canonical metric, with possible implications for the study of stability.
We also show that on a Fano K\"ahler-Einstein manifold all of the functionals $E_k$ are bounded below on the space of metrics with nonnegative Ricci curvature.

\bigskip
\noindent {\it Mathematical Subject Classification: 32Q20 (Primary),
53C21 (Secondary)} \setlength\arraycolsep{2pt}
\addtocounter{section}{1}

\bigskip
\bigskip
\noindent
{\bf 1. Introduction}

\bigskip

The problem of finding necessary and sufficient conditions for the existence of extremal metrics, which include K\"ahler-Einstein metrics, on a compact K\"ahler manifold $M$ has been the subject of intense study over the last few decades and is still largely open.  If $M$ has zero or negative first Chern class then it is known by the work of Yau \cite{Ya1} and Yau, Aubin \cite{Ya1}, \cite{Au} that $M$ has a K\"ahler-Einstein metric.  When $c_1(M)>0$, so that $M$ is Fano, there is a well-known conjecture of Yau \cite{Ya2} that the manifold admits a K\"ahler-Einstein metric if and only if it is stable in the sense of geometric invariant theory.

There are now several different notions of stability for manifolds \cite{Ti2}, \cite{PhSt1}, \cite{Do2}, \cite{RoTh}.  Donaldson showed that the existence of a constant scalar curvature metric is sufficient for the manifold to be asymptotically Chow stable \cite{Do1} (under an assumption on the automorphism group).
It is conjectured by Tian \cite{Ti2} that the existence of a K\"ahler-Einstein metric should be equivalent to his `K-stability'.   This stability is defined in terms of the Futaki invariant \cite{Fu}, \cite{DiTi} of the central fiber of degenerations of the manifold. Donaldson  \cite{Do2}  introduced a variant of K-stability extending Tian's definition.

The behavior of Mabuchi's \cite{Mb} energy functional  is central to this problem.  It was shown by Bando and Mabuchi \cite{BaMb}, \cite{Ba} that if a Fano manifold admits a K\"ahler-Einstein metric then the Mabuchi energy is bounded below. Recently, it has been shown by Chen and Tian \cite{ChTi3} that if $M$ admits an extremal metric in a given class then the (modified) Mabuchi energy is bounded below in that class.  Donaldson has given an alternative proof for constant scalar curvature metrics with a condition on the space of automorphisms \cite{Do3}.  Moreover, if a lower bound on the  Mabuchi energy is given then the class is K-semistable \cite{Ti2}, \cite{PaTi}.  Conversely, Donaldson \cite{Do2} showed that, for toric surfaces, K-stability implies the lower boundedness of the Mabuchi energy.

In addition, the existence of a K\"ahler-Einstein metric on a Fano manifold has been shown to be equivalent to the `properness' of the Mabuchi energy \cite{Ti2}.  Tian conjectured \cite{Ti3} that the existence of a constant scalar curvature K\"ahler metric be equivalent to this condition on the Mabuchi energy.  This holds when the first Chern class is a multiple of the K\"ahler class.  (If $c_1(M)<0$, it has been shown  \cite{Ch1}, \cite{We}, \cite{SoWe} that the Mabuchi energy is proper on certain classes which are not multiples of the canonical class.  It is not yet known whether there exists a constant scalar curvature metric in such classes.)

In this paper, we discuss a family of functionals $E_k$, for $k=0, \ldots, n$,  which were introduced by Chen and Tian \cite{ChTi1}.  They are generalizations of the Mabuchi energy, with $E_0$ being precisely Mabuchi's functional.  The functionals $E_k$ can be described in terms of the Deligne pairing \cite{De}.  This construction of Deligne has provided a very useful way to understand questions of stability \cite{Zh}, \cite{PhSt1}, \cite{PhSt2}.  Phong and Sturm \cite{PhSt3} show that, up to a normalization term, the Mabuchi energy corresponds to the Deligne pairing $\langle L, \ldots, L, K^{-1} \rangle$.
Generalizing this, the functionals $E_k$ can be described in terms of the pairing:
$$\langle \stackrel{n-k}{\overbrace{L, \ldots, L}}, \stackrel{k+1}{\overbrace{K^{-1}, \ldots, K^{-1}}} \rangle.$$
The fact that the functionals $E_k$ can be formulated in this way seems now to be known by some experts in the field, and was pointed out to us by Jacob Sturm in 2002.  However, since it does not appear in the literature, we have included a short explanation of this correspondence (see section 2).

A critical metric $\omega$ of $E_k$ is a solution of the equation
\begin{eqnarray*}
\sigma_{k+1} (\omega) - \Delta (\sigma_k(\omega))  =  \textrm{constant},
\end{eqnarray*}
where  $\sigma_k(\omega)$ is the $k$th elementary symmetric polynomial in the eigenvalues of the Ricci tensor of the metric $\omega$.  Critical points for $E_0$ are precisely the constant scalar curvature metrics. Critical
points for $E_n$ are the metrics of constant central curvature as described
by
 Maschler \cite{Ms}.  K\"ahler-Einstein metrics are solutions to the above equation for all $k$.  The critical metrics are discussed more in section 2.

The functionals $E_k$ were used by Chen and Tian \cite{ChTi1, ChTi2} to obtain convergence of the normalized K\"ahler-Ricci flow on K\"ahler-Einstein manifolds with positive bisectional curvature (see \cite{PhSt4} for a related result).  The Mabuchi energy is decreasing along the flow.  The functional $E_1$ is also decreasing, as long as the sum of the Ricci curvature and the metric is nonnegative.  A major part of the argument in \cite{ChTi1} is to show that the $E_k$ can be bounded from below along the K\"ahler-Ricci flow assuming nonnegative Ricci curvature along the flow and the existence of a K\"ahler-Einstein metric.

In a recent preprint, Chen \cite{Ch2} has proved a stability result for $E_1$ for Fano manifolds in the sense of the K\"ahler-Ricci flow\footnote{Shortly after this paper was first posted we learned that, in an unpublished work \cite{Ch3}, Chen has proved the following: if there exists a K\"ahler-Einstein metric then $E_1$ is bounded below along the K\"ahler-Ricci flow.}.   In \cite{Ch2}, Chen asked whether $E_1$ is bounded below or proper on the full space of potentials (not just along the flow) if there exists a K\"ahler-Einstein metric.  In this paper we answer Chen's question:  $E_1$ is bounded below if there exists a K\"ahler-Einstein metric, and the lower bound is attained by this metric.  Moreover, modulo holomorphic vector fields, $E_1$ is proper if and only if there exists a K\"ahler-Einstein metric.  We also show that, again assuming the existence of a K\"ahler-Einstein metric, the functionals $E_k$  are bounded below on the space of metrics with nonnegative Ricci c
 urvature.

We now state these results more precisely. Let $\omega$ be a K\"ahler form on the compact manifold $M$ of complex dimension $n$.  Write $P(M, \omega)$ for the space of all smooth functions $\phi$ on $M$ such that
$$\omega_{\phi} = \omega + \ddbp >0.$$
For $\phi$ in $P(M, \omega)$, let $\phi_t$ be a path in $P(M, \omega)$ with $\phi_0 = 0$ and $\phi_1= \phi$.  The functional $E_{k, \omega}$ for $k=0, \ldots, n$ is defined by
\begin{eqnarray} \nonumber
E_{k,\omega} (\phi) & = & \frac{k+1}{V} \int_0^1 \int_{M} (\Delta_{\phi_t} \dotp_t ) \, \Rc (\omega_{\phi_t})^k \wedge \omega_{\phi_t}^{n-k}dt \\ \label{eqnEk}
&& \hspace{-15pt} \mbox{} - \frac{n-k}{V} \int_0^1 \int_M \dotp_t \, ( \Rc (\omega_{\phi_t})^{k+1} - \mu_k \, \omega_{\phi_t}^{k+1}) \wedge \omega_{\phi_t}^{n-k-1}dt,
\end{eqnarray}
where $V$ is the volume $\int_M \omega^n$, and $\mu_k$ is the constant, depending only on the classes $[\omega]$ and $c_1(M)$ given by
$$\mu_k = \frac{\int_M \Rc(\omega)^{k+1} \wedge \omega^{n-k-1}}{\int_M \omega^n} = (2\pi)^{k+1}\frac{ [ K^{-1}]^{k+1} \cdot [\omega]^{n-k-1}}{[\omega]^n}.$$
The functional is independent of the choice of path.
We will often write $E_k(\omega, \omega_{\phi})$ instead of $E_{k, \omega} (\phi)$.

In fact, Chen and Tian \cite{ChTi1} first define $E_k$ by a different (and explicit) formula, making use of  a  generalization of the Liouville energy, which they call $E^0_k$.  The Deligne pairing provides another explicit formula (Proposition \ref{propdeligne}).

Suppose now that $M$ has positive first Chern class and denote by $\kahler$ the space of K\"ahler metrics in $2\pi c_1(M)$.  Notice that for $\omega$ in $\kahler$, the corresponding constant $\mu_k$ is equal to 1.
We have the following theorem on the lower boundedness of the functionals $E_k$.

\begin{theorem} \label{theorem1} Let $(M, \omega_{KE})$ be a K\"ahler-Einstein manifold with $c_1(M)>0$.  Then, for $k=0, \ldots, n$, and for all $\tilde{\omega} \in \kahler$ with $\emph{Ric}(\tilde{\omega}) \ge 0$,
$$E_{k}(\omega_{KE}, \tilde{\omega})  \ge 0,$$
and equality is attained if and only if $\tilde{\omega}$ is a K\"ahler-Einstein metric.
\end{theorem}

In the case of $E_1$ we obtain lower boundedness on the whole space $\kahler$.  In addition, it is an easy result that for a Calabi-Yau manifold, $E_1$ is bounded below on every class.  Putting these two cases together we obtain:

\pagebreak[3]
\begin{theorem} \label{theorem2} Let $(M, \omega_{KE})$ be a K\"ahler-Einstein manifold with $c_1(M) > 0$ or $c_1(M)=0$.  Then for all K\"ahler metrics $\omega'$ in the class $[\omega_{KE}]$,
$$E_1(\omega_{KE}, \omega') \ge 0,$$
and equality is attained if and only if $\omega'$ is a K\"ahler-Einstein metric.
\end{theorem}

We show that if $(M, \omega_{KE})$ is K\"ahler-Einstein with $c_1(M)>0$ and if there are no holomorphic vector fields, $E_1$ is bounded below by the Aubin-Yau energy functional $J$ raised to a small power.  This implies that $E_1$ is proper on $P(M, \omega_{KE})$.  If there exist holomorphic vector fields, then the statement changes slightly (c.f. \cite{Ti2}).

\begin{theorem} \label{theorem3}
Let $(M, \omega_{KE})$ be a compact K\"ahler-Einstein manifold with $c_1(M) >0$.   Then there exists $\delta$ depending only on $n$ such that the following hold:
\begin{enumerate}
\item[(i)] If $M$ admits no nontrivial holomorphic vector fields then there exist positive constants $C$ and $C'$ depending only on $\omega_{KE}$ such that for all $\theta$ in $P(M, \omega_{KE})$,
$$E_{1,\omega_{KE}} (\theta) \ge C J_{\omega_{KE}} (\theta)^{\delta} - C'.$$
\item[(ii)] In general, let $G$ be the maximal compact subgroup of $\emph{Aut}^0(M)$ which fixes $\omega_{KE}$, where $\emph{Aut}^0(M)$ is
 the component of the automorphism group containing the identity.  Then there exist positive constants $C$ and $C'$ depending only on $\omega_{KE}$ such that for all $\theta$ in $P_G(M, \omega_{KE})$,
$$E_{1,\omega_{KE}} (\theta) \ge C J_{\omega_{KE}} (\theta)^{\delta} - C',$$
where $P_G(M, \omega_{KE})$ consists of the $G$-invariant elements of $P(M, \omega_{KE})$.
\end{enumerate}
\end{theorem}

\begin{remark}
The analagous result above is proved in \cite{Ti2}, \cite{TiZh} for the $F$ functional \cite{Di}, giving a generalized Moser-Trudinger inequality.  We use a similar argument.  The corresponding inequality for the Mabuchi energy also holds \cite{Ti3} since, up to a constant, it can be bounded below by $F$.  It would be interesting to find the best constant $\delta=\delta(n)$ for which these inequalities hold.
With some work, modifying the argument in \cite{Ti2}, one can show that $\delta$ can be taken to be arbitrarily close to $1/(4n+1)$, but we doubt that this is optimal.
\end{remark}

\begin{remark}  We hope that the above results on the $E_1$ functional may have some applications to the study of the stability of $M$.  Indeed,
let $\pi_1: \mathcal{X}\rightarrow Z$ be an $SL(N+1,
\mathbf{C})$-equivariant holomorphic fibration between smooth
varieties such that $\mathcal{X}\subset Z\times\mathbf{CP}^N$ is a
family of subvarieties of dimension $n$ with an action
of $SL(N+1, \mathbf{C})$ on $\mathbf{CP}^N$. Tian defines
CM-stability \cite{Ti2} for $X_z=\pi_1^{-1}(z)$ in terms of the virtual bundle:
$$\mathcal{E}=(\mathcal{K}^{-1}-\mathcal{K})
\otimes(L-L^{-1})^n-\frac{n\mu_0}{n+1}(L-L^{-1})^{n+1},$$ where
$\mathcal{K}=K_{\mathcal{X}}\otimes K_Z^{-1}$ is the relative
canonical bundle and $L$ is  the pullback of the hyperplane bundle on
$\mathbf{CP}^N$ via the second projection $\pi_2$ (alternatively, one can use the language of the Deligne pairing).
When $X_z$ is Fano,  Tian proved that $X_z$ is weakly CM-stable if
it is K\"ahler-Einstein, using the properness of the Mabuchi
energy $E_0$.   One can define a similar notion of stability for $X_z$ with
respect to the virtual bundle
$$\mathcal{E}_k=(\mathcal{K}^{-1}-\mathcal{K})^{k+1}
\otimes(L-L^{-1})^{n-k}-\frac{(n-k)\mu_k}{n+1}
(L-L^{-1})^{n+1}.$$ For $k=1$, one might guess that
$X_z$ is stable if it is K\"ahler-Einstein, since $E_1$ is proper. It would be interesting to try to relate this notion of stability to an analogue of K-stability expressed in terms of the holomorphic invariants $\mathcal{F}_k$ \cite{ChTi1} which generalize the Futaki invariant.
\end{remark}

We also have a converse to Theorem \ref{theorem3}.

\begin{theorem} \label{theorem4}
Let $(M, \omega)$ be a compact K\"ahler manifold with $c_1(M)>0$.  Suppose that $\omega \in 2\pi c_1(M)$. Then the following hold:
\begin{enumerate}
\item[(i)] Suppose that $(M, \omega)$ admits no nontrivial holomorphic vector fields.  Then $M$ admits a K\"ahler-Einstein metric if and only if $E_1$ is proper on $P(M, \omega)$.
\item[(ii)] In general, suppose that $\omega$ is invariant under $G$, a maximal compact subgroup of $\emph{Aut}^0(M)$.  Then $M$ admits a $G$-invariant K\"ahler-Einstein metric if and only if $E_1$ is proper on $P_{G} (M, \omega)$.\end{enumerate}
\end{theorem}

This gives  a new analytic condition for a Fano manifold to admit a K\"ahler-Einstein metric.  Indeed, together with the result of Tian \cite{Ti2}, at least modulo holomorphic vector fields, we have:
$$
\begin{array}{ccccc}
&& $M$ \textrm{ admits a} && \\
E_1 \textrm{ is proper } & \quad \Longleftrightarrow \quad & \textrm{K\"ahler-Einstein}  &  \quad \Longleftrightarrow \quad & \textrm{ Mabuchi energy is proper} \\
&& \textrm{metric}&&
\end{array}$$
and one might expect some versions of stability to be equivalent to these as well.

\begin{remark}
It is natural to ask whether there exist critical metrics for $E_k$ which are not K\"ahler-Einstein.  Chen and Tian \cite{ChTi1} observed that for $k=n$ the only critical metrics with positive Ricci curvature are K\"ahler-Einstein.
Maschler \cite{Ms} proved this result without the assumption on the Ricci curvature when $c_1(M)>0$ or $c_1(M)<0$.  We see from Theorem \ref{theorem2} that on a K\"ahler-Einstein manifold with $c_1(M)>0$, a critical metric for $E_1$ which is not K\"ahler-Einstein could not give an absolute minimum of $E_1$.
\end{remark}

In section 2, we describe some of the properties of the $E_k$ functionals and the Aubin-Yau functionals $I$ and $J$.  In sections 3 and 4 we prove Theorems \ref{theorem1} and \ref{theorem2} respectively.  Our general method follows that of Bando and Mabuchi \cite{BaMb} and Bando \cite{Ba}.  They consider the Monge-Amp\`ere equations (\ref{eqnma1}) and (\ref{eqnma2}) respectively.  The Mabuchi energy is decreasing in $t$ for both of these equations, and one can solve them backwards to show that the Mabuchi energy is bounded below on the space of metrics with positive Ricci curvature \cite{BaMb}, and on all metrics in the given class \cite{Ba}.  The functionals $E_k$ may not be decreasing in $t$.  However, we are able to show, with some extra calculation, that the value at $t=0$ is greater than the value at $t=1$.  This works for all the $E_k$ for the first equation, but only seems to work for $E_1$ for the second.

Finally, in section 5 we adapt the results of Tian \cite{Ti2} and Tian-Zhu \cite{TiZh}, making use of our formulas from sections 3 and 4, to prove Theorem \ref{theorem3} and Theorem \ref{theorem4}.

\setcounter{equation}{0}
\addtocounter{section}{1}
\bigskip
\bigskip
\pagebreak[3]
\noindent
{\bf 2. Properties of the $E_k$ and $I$, $J$ functionals}
\bigskip

\noindent
{\it 2.1  Basic properties of $E_k$}
\bigskip

\noindent
It is immediate from the formula (\ref{eqnEk}) that the functionals $E_k$ satisfy the cocycle condition
$$E_k (\omega_1, \omega_2) + E_k (\omega_2, \omega_3) + E_k (\omega_3, \omega_1) = 0,$$
for any metrics $\omega_1$, $\omega_2$ and $\omega_3$ in the same K\"ahler class.  We will make use of this relation in sections 3 and 4.  Observe also that $E_{k, \omega}$ is invariant under addition of constants:
$$E_{k,\omega} (\phi + c) = E_{k, \omega}(\phi)$$
for all $\phi$ in $P(M, \omega)$ and any constant $c$.

\bigskip
\pagebreak[3]
\noindent
{\it 2.2  Critical metrics}
\bigskip

\noindent
The $E_k$ functionals can be written as
\begin{eqnarray*}
E_{k, \omega} (\phi)
&  = & \frac{1}{V} \int_0^1 \int_M \dot{\phi}_t \left[ (k+1) \binom{n}{k}^{-1} (\Delta_{\phi_t} (\sigma_k(\omega_{\phi_t})) - \sigma_{k+1}(\omega_{\phi_t})) + (n-k)\mu_k \right] \omega_{\phi_t}^n dt,
\end{eqnarray*}
where $\sigma_k (\omega)$ is given by
$$(\omega + t\,  \Rc(\omega))^n = \left( \sum_{k=0}^n \sigma_k(\omega) \, t^k \right) \omega^n.$$
If at a point $p$ on $M$ we pick a normal coordinate system  in which the Ricci tensor is diagonal with entries $\lambda_1, \ldots, \lambda_n$ then at $p$, we have $\sigma_{k}(\omega) = \sum_{i_1 < \cdots < i_k} \lambda_{i_1} \lambda_{i_2} \cdots \lambda_{i_k}.$
For example, $\sigma_0=1$, $\sigma_1=R$ and $\sigma_2 = \frac{1}{2} (R^2 - |\Rc|^2)$.
The critical points of $E_k$ are metrics $\omega$ satisfying
\begin{equation} \label{eqncritical}
\sigma_{k+1} (\omega) - \Delta (\sigma_k(\omega))  =  \binom{n}{k+1} \mu_k,
\end{equation}
The critical points for $k=0$ are constant scalar curvature metrics. For $k=n$, the critical points satisfy $(\textrm{Ric}(\omega))^n = c\,  \omega^n$
for some constant $c$.  These metrics are said to have constant central curvature
and are studied in depth in \cite{Ms}, \cite{HwMs}.


\bigskip
\noindent
{\it 2.3  Deligne Pairing}
\bigskip

\noindent
We will now describe the functionals $E_k$ in terms of the Deligne pairing.
Let $\pi: \mathcal{X} \rightarrow S$ be a flat projective morphism of integral schemes of relative dimension $n$.  For each $s$, $X_s = \pi^{-1}(s)$ is a projective variety of dimension $n$.  Let $\mathcal{L}_0, \ldots, \mathcal{L}_n$ be Hermitian line bundles on $\mathcal{X}$.  We denote the corresponding Deligne pairing - a Hermitian line bundle on $S$ - by
$$\langle \mathcal{L}_0, \ldots, \mathcal{L}_n \rangle (\mathcal{X}/S).$$
For its definition, we refer the reader to the references \cite{De}, \cite{Zh} and \cite{PhSt3}.
For a smooth function $\phi$ on $\mathcal{X}$ (or $S$) denote by $\mathcal{O}(\phi)$ the trivial line bundle on $\mathcal{X}$ (or $S$) with metric $e^{-\phi}$.  The only property of the Deligne pairing we will use is the change of metric formula:
\begin{equation} \nonumber
\langle \mathcal{L}_0 \otimes \mathcal{O}(\phi_0), \ldots, \mathcal{L}_n \otimes \mathcal{O} (\phi_n) \rangle (\mathcal{X}/S) = \langle \mathcal{L}_0, \ldots, \mathcal{L}_n \rangle (\mathcal{X}/S) \otimes \mathcal{O}(E), \label{eqnchange}
\end{equation}
where $E$ is the function on $S$ given by
$$E(s) = \int_{X_s} \sum_{j=0}^n \phi_j \, \left( \prod_{i<j} c'_1 (\mathcal{L}_i \otimes \mathcal{O}(\phi_i)) \right) \wedge \left( \prod_{i>j} c'_1 (\mathcal{L}_i) \right),$$
for $c'_1(\mathcal{L}) = - \frac{1}{2\pi} \ddbar \log \| \cdot \|_{\mathcal{L}}$.

In our case, $\mathcal{X}$ will be the variety $M$ and $S$ will be a single point. In the following we will omit any reference to $\mathcal{X}$ or $S$. We fix a metric $\omega \in 2\pi c_1(L)$ for some line bundle $L$ on $M$, and an element $\phi$ of $P(M, \omega)$.  Let $\mathcal{L}$ be the Hermitian line bundle $(L, h)$ with
$-\ddbar \log h = \omega$.  Define $\tilde{\mathcal{L}}$ by $\tilde{\mathcal{L}}= \mathcal{L} \otimes \mathcal{O}(\phi).$  Let $\AK$ be the anticanonical bundle $K^{-1}$ equipped with the metric $\omega^n$, and let $\widetilde{\AK}$ be the same bundle with the metric $\omega_{\phi}^n$.  In other words,
$$\widetilde{\AK} = \AK \otimes \mathcal{O} (\log \frac{\omega^n}{\omega_{\phi}^n}).$$
Then we have the following formula for the $E_k$ functional.

\begin{proposition} \label{propdeligne}
\begin{eqnarray*}
\lefteqn{\langle \stackrel{n-k}{\overbrace{\tilde{\mathcal{L}}, \ldots, \tilde{\mathcal{L}}}}, \stackrel{k+1}{\overbrace{\widetilde{\AK}, \ldots, \widetilde{\AK}}} \rangle \otimes ( \langle \tilde{\mathcal{L}}, \ldots, \tilde{\mathcal{L}} \rangle )^{- \frac{(n-k) \mu_k}{n+1}} } \\
&& = \langle \stackrel{n-k}{\overbrace{\mathcal{L}, \ldots, \mathcal{L}}}, \stackrel{k+1}{\overbrace{\AK, \ldots, \AK}} \rangle \otimes (\langle \mathcal{L}, \ldots, \mathcal{L} \rangle )^{- \frac{(n-k) \mu_k}{n+1}} \otimes \mathcal{O}(- \frac{V}{(2\pi)^n} E_{k, \omega}(\phi) )
\end{eqnarray*}

\end{proposition}

\begin{proof}  Define a functional $a_{k, \omega}$ by
\begin{equation} \label{eqnak}
\langle \stackrel{n-k}{\overbrace{\tilde{\mathcal{L}}, \ldots, \tilde{\mathcal{L}}}}, \stackrel{k+1}{\overbrace{\widetilde{\AK}, \ldots, \widetilde{\AK}}} \rangle
 =  \langle \stackrel{n-k}{\overbrace{\mathcal{L}, \ldots, \mathcal{L}}}, \stackrel{k+1}{\overbrace{\AK, \ldots, \AK}} \rangle \otimes \mathcal{O} \left( \frac{1}{(2\pi)^n} a_{k, \omega}(\phi) \right).
\end{equation}
From the change in metric formula we see that
\begin{eqnarray} \nonumber
a_{k,\omega} (\phi) & = &  \sum_{j=0}^{n-k-1} \int_M \phi \, \omega_{\phi}^j \wedge \Rc(\omega)^{k+1} \wedge \omega^{n-j-k-1} \\ \nonumber
&& \mbox{} + \sum_{j=0}^k \int_M \log \left(\frac{\omega^n}{\op^n} \right) \Rc(\op)^j \wedge \op^{n-k} \wedge \Rc(\omega)^{k-j}.
\end{eqnarray}
Now calculate
\allowdisplaybreaks{
\begin{eqnarray}  \nonumber
\frac{d}{dt} a_{k,\omega} (\phi)
& = & \sum_{j=0}^{n-k-1}\int_M \dot{\phi} \, \op^j \wedge \Rc(\omega)^{k+1} \wedge \omega^{n-j-k-1} \\  \nonumber
&& \mbox{} + \sum_{j=0}^{n-k-1} j \int_M \phi \, \op^{j-1} \wedge \ddbar \dot{\phi} \wedge \Rc(\omega)^{k+1} \wedge \omega^{n-j-k-1} \\ \nonumber
&& \mbox{} - \sum_{j=0}^k \int_M (\Delta_{\phi} \dot{\phi}) \Rc(\op)^j \wedge \op^{n-k} \wedge \Rc(\omega)^{k-j} \\ \nonumber
&& \mbox{} - \sum_{j=0}^k j \int_M \log \left( \frac{\omega^n}{\op^n} \right) \Rc(\op)^{j-1} \wedge \ddbar (\Delta_{\phi} \dot{\phi}) \wedge \op^{n-k} \wedge \Rc(\omega)^{k-j} \\ \nonumber
&& \mbox{} + \sum_{j=0}^k (n-k) \int_M \log \left( \frac{\omega^n}{\op^n} \right) \Rc(\op)^j \wedge \op^{n-k-1} \wedge \ddbar \dot{\phi} \wedge \Rc(\omega)^{k-j}.
\end{eqnarray}}
Integrating by parts, and rearranging, we have
\begin{eqnarray} \nonumber
\frac{d}{dt} a_{k,\omega} (\phi)
& = &  -(k+1) \int_M (\Delta_{\phi} \dot{\phi}) \Rc(\op)^{k} \wedge \op^{n-k} \\ \label{eqnak2}
&& \mbox{} + (n-k) \int_M \dot{\phi} \, \Rc(\op)^{k+1} \wedge \op^{n-k-1}.
\end{eqnarray}
We recognize this (after dividing by $-V$) as the first two terms of the derivative of $E_{k, \omega}$.  We are left with the term involving the constant $\mu_k$.   As in \cite{Zh}, define a functional $b_{\omega}$ by
\begin{equation} \label{eqnb}
\langle \tilde{\mathcal{L}}, \ldots, \tilde{\mathcal{L}} \rangle = \langle \mathcal{L}, \ldots, \mathcal{L} \rangle \otimes \mathcal{O}\left( \frac{1}{(2\pi)^n} b_{\omega} (\phi) \right).
\end{equation}
From the change in metric formula, we see that $b_{\omega}$ is the well-known functional:
$$b_{\omega}(\phi) =  \sum_{i=0}^n \int_M \phi \, \op^i \wedge \omega^{n-i}.$$
Its derivative along a path in $P(M, \omega)$ is given by
\begin{eqnarray} \label{eqnb2}
\frac{d}{dt} b_{\omega} (\phi) & =&  (n+1) \int_M \dot{\phi} \, \op^n.
\end{eqnarray}
Comparing with (\ref{eqnEk}), we see from (\ref{eqnak2}) and (\ref{eqnb2}) that
$$E_{k, \omega} (\phi) = - \frac{1}{V} a_{k,{\omega}}(\phi) + \frac{(n-k) \mu_k} {V(n+1)} b_{\omega} (\phi).$$
The proposition then follows from this along with the defining equations (\ref{eqnak}) and (\ref{eqnb}). Q.E.D.\end{proof}

\bigskip
\pagebreak[3]
\noindent
{\it 2.4 The $I$ and $J$ functionals}
\bigskip

The Aubin-Yau energy functionals $I$ and $J$ are defined as follows.  For $\phi \in P(M,\omega)$, set
\begin{eqnarray*}
I_{\omega} (\phi) & = & \frac{1}{V} \sum_{i=0}^{n-1} \int_M \sqrt{-1} \partial \phi \wedge \dbar \phi \wedge \omega^i \wedge \omega_{\phi}^{n-1-i}, \\
J_{\omega} (\phi) & = & \frac{1}{V} \sum_{i=0}^{n-1} \frac{i+1}{n+1} \int_M \sqrt{-1} \partial \phi \wedge \dbar \phi \wedge \omega^i \wedge \omega_{\phi}^{n-1-i}.
\end{eqnarray*}
Their difference is given by
$$
(I_{\omega} - J_{\omega})(\phi) =  \frac{1}{V} \sum_{i=0}^{n-1} \frac{n-i}{n+1} \int_M \sqrt{-1} \partial \phi \wedge \dbar \phi \wedge \omega^i \wedge \omega_{\phi}^{n-1-i}.
$$
Observe that $I$, $J$ and $I-J$ are all nonnegative and equivalent.
A well-known calculation shows that for any family of functions $\phi=\phi_t \in P(M, \omega)$,
\begin{equation} \label{eqnddtIminusJ}
\frac{d}{dt} (I_{\omega} - J_{\omega})(\phi_t) = - \frac{1}{V} \int_M \phi_t (\Delta_{\phi_t} \dot{\phi_t}) \omega_{\phi_t}^n.
\end{equation}

We now define the notion of properness for a functional.
Following Tian \cite{Ti3}, we say that a functional $T$ on a $P(M, \omega)$ is \emph{proper} if there
there exists an
increasing function
$f: [0, \infty) \rightarrow \mathbf{R},$
satisfying $f(x) \rightarrow \infty$ as $x \rightarrow \infty$, such that for any $\phi \in P(M,\omega)$,
$T(\phi) \ge f\left( J_{\omega} (\phi) \right).$  The condition of properness is independent of the choice of metric $\omega$ in the class.

\addtocounter{section}{1}
\setcounter{equation}{0}
\bigskip
\bigskip
\pagebreak[3]
\noindent
{\bf 3.  The lower boundedness of the $E_k$ functionals}
\bigskip

\noindent
In the following, we will denote by $\kahlerp$ the set of metrics $\omega$ in $\kahler$ with strictly positive Ricci curvature.  We will write $\mathcal{E}$ for the space of metrics $\omega \in \kahler$ with $\Rc(\omega) = \omega$.
We first show that the functionals $E_k$ are bounded below on $\kahlerp$.   We need the following result from \cite{BaMb} (see also \cite{Si} for a good exposition of their work):

\begin{lemma} \label{lemmaBM}
Let $\tilde{\omega} \in \kahlerp$ and set $\omega_0 = \emph{Ric}(\tilde{\omega})$.
If $\mathcal{E}$ is non-empty, there exists a K\"ahler-Einstein metric $\omega_{KE}' \in \mathcal{E}$ and, for $\delta>0$ sufficiently small, a family of metrics $\omega_0^{\epsilon} \in \kahlerp$ for $\ep \in (0,\delta)$ such that:
\begin{enumerate}
\item[(i)] $\omega_0^{\ep} \rightarrow \omega_0$ in $C^{\infty}$ as $\epsilon \rightarrow 0$.
\item[(ii)] Let $f^{\ep}$ be the unique smooth function on $M$ satisfying
$$\emph{Ric}(\omega_0^{\epsilon}) - \omega_0^{\epsilon} = \ddbar f^{\epsilon} \quad \textrm{and} \quad \int_M e^{f^{\epsilon}} (\omega_0^{\epsilon})^n = \int_M  (\omega_0^{\epsilon})^n.$$
Then for each $\epsilon \in (0, \delta)$, there exists a smooth solution $\psi^{\epsilon}_t$ (for $0 \le t \le 1$) of the equation
$$ (\omega_0^{\epsilon} + \ddbar \psi^{\epsilon}_t)^n = e^{-t \psi^{\epsilon}_t + f^{\ep}} (\omega_0^{\epsilon})^n,$$
with
$\omega_0^{\epsilon} + \ddbar \psi^{\epsilon}_1 = \omega_{KE}'.$
\item[(iii)]  Set $\tilde{\omega}^{\epsilon} = \omega_0^{\epsilon} + \ddbar \psi^{\epsilon}_0$.  Then
$\tilde{\omega}^{\epsilon} \rightarrow \tilde{\omega}$ in  $C^{\infty}$ as  $\epsilon \rightarrow 0.$

\end{enumerate}
\end{lemma}

In order to make use of this for the proof of Theorem \ref{theorem1}, we need the following lemma:

\begin{lemma} \label{lemmaEk}
Let $\omega$ be a metric in $\kahler$, and let $\phi_t$ be a solution of
\begin{equation} \label{eqnma1}
\omega_{\phi_t}^n = e^{-t\phi_t + f} \omega^n,
\end{equation}
where $f$ is given by
$$\emph{Ric}(\omega) - \omega = \ddbar f \quad \textrm{and} \quad \int_M e^{f} \omega^n = V.$$
Then $E_k (\omega, \omega_{\phi_0})  \ge E_k(\omega, \omega_{\phi_1}).$
\end{lemma}

Before proving this lemma, we will complete the first part of the
proof of Theorem \ref{theorem1} in the case $\tilde{\omega} \in
\kahlerp$. From the above lemmas, for any such $\tilde{\omega}$,
$$E_{k} (\omega_0^{\epsilon}, \tilde{\omega}^{\epsilon}) \ge E_k(\omega_0^{\epsilon}, \omega_{KE}').$$
By the cocycle condition,
$$E_{k} (\omega_{KE}, \tilde{\omega}^{\epsilon}) \ge E_k(\omega_{KE}, \omega_{KE}').$$
Letting $\epsilon \rightarrow 0$, we obtain from the continuity of
$E_k$,
$$E_{k} (\omega_{KE}, \tilde{\omega}) \ge E_k(\omega_{KE}, \omega_{KE}').$$
We need a lemma to show that the right hand side of this equation is
zero.

\begin{lemma} \label{lemmainvariants}
 Let $\omega_{KE}$ and $\omega'_{KE}$ be  in  $\mathcal{E}$.  Then
$E_k(\omega_{KE}, \omega'_{KE}) = 0.$
\end{lemma}
\begin{proof}
Chen and Tian \cite{ChTi1} define an invariant of a Fano manifold $M$ which generalizes the Futaki invariant.  For each holomorphic vector field $X$, they define for any $\omega \in \kahler$,
\begin{eqnarray*}
\lefteqn{\mathcal{F}_k (X)= (n-k) \int_M h_X \omega^n } \\
&& \mbox{} + \int_M \left( (k+1) \Delta_{\omega} h_X \Rc(\omega)^k \wedge \omega^{n-k} - (n-k) h_X \Rc(\omega)^{k+1} \wedge \omega^{n-k-1}\right),
\end{eqnarray*}
where $h_X$ is a function (unique up to constants) satisfying $\mathcal{L}_X (\omega) = \ddbar h_X.$
They prove that $\mathcal{F}_k (X)$ does not depend on the choice of $\omega$ in $\kahler$, and that if $\mathcal{E}$ is non-empty, the invariant vanishes.  Moreover, it is shown that if $X$ is a holomorphic vector field and $\{ \Phi(t) \}$ the one-parameter subgroup of automorphisms induced by $\textrm{Re}(X)$, then for any $\omega$,
$$\frac{d}{dt} E_k (\omega, \Phi_t^* \omega) = \frac{1}{V} \textrm{Re} (\mathcal{F}_k(X)).$$
Given this, the lemma is an immediate consequence of the theorem of Bando and Mabuchi \cite{BaMb} that the space $\mathcal{E}$ is a single orbit of the action of the group of holomorphic automorphisms of $M$.    Q.E.D.
\end{proof}

We will now prove Lemma \ref{lemmaEk}.

\bigskip
\noindent
{\bf Proof of Lemma \ref{lemmaEk}}
\bigskip

Let $\phi_t$ be a solution of  (\ref{eqnma1}).
Differentiating (\ref{eqnma1}) with respect to $t$, we have
\begin{equation} \label{eqnlapphit}
\Delta_{\phi_t} \dotp_t = -t\dotp_t -\phi_t.
\end{equation}
Applying the operator $- \ddbar \log$ to (\ref{eqnma1}) we obtain
$$\Rc(\omega_{\phi_t}) = \Rc(\omega) + t \ddbp_t - \ddbar f,$$
which, from the definition of $f$, can be written
\begin{equation} \label{eqnric}
\Rc (\omega_{\phi_t}) = \omega_{\phi_t} + (t-1) \ddbp_t.
\end{equation}

We will need to make use of the following fact which is well-known:  for $\phi_t$ a solution of (\ref{eqnma1}), the first eigenvalue $\lambda_1$ of $\Delta_{\phi_t}$ satisfies $\lambda_1 \ge t$ if $0 \le t  \le 1$ (and the inequality is strict if $t<1$.)  An immediate consequence of this and (\ref{eqnlapphit}) is the following inequality:
\begin{equation} \label{eqnintegralinequality}
\int_M \phi_t (\Delta_{\phi_t} \dot{\phi}_t) \omega_{\phi_t}^n \le 0,
\end{equation}
which by (\ref{eqnddtIminusJ}) is equivalent to the fact that $(I_{\omega}-J_{\omega})(\phi_t)$ is increasing in $t$.

We will find a formula for the expression $(E_{k, \omega} (\phi_1) - E_{k,\omega} (\phi_0))$ by calculating $\int_0^1 \frac{d}{dt} E_{k, \omega} (\phi_t) dt.$
Lemma \ref{lemmaEk} and the first part of Theorem \ref{theorem1} will then follow from (\ref{eqnintegralinequality}) and the next lemma, which is the key result of this section. In the statement and proof of the lemma, we will simplify the notation by omitting the subscript $t$.

\begin{lemma} \label{lemmaEkformula1}
Let $\phi= \phi_t$ be a solution of (\ref{eqnma1}).  Then
\begin{eqnarray*}
E_{k, \omega} (\phi_1) - E_{k, \omega} (\phi_0) & = &
\frac{k+1}{V} \int_0^1 \int_M (1-t) \phi (\Delta_{\phi} \dotp) \op^n dt \\
&& \left. \mbox{} -  \frac{1}{V} \sum_{i=0}^{k-1} (k-i) \int_M \sqrt{-1} \partial \phi \wedge \dbar \phi \wedge\omega^i \wedge \op^{n-i-1} \right|_{t=0}.\\
\end{eqnarray*}
\end{lemma}
\begin{proof} Using (\ref{eqnlapphit}) and (\ref{eqnric}), calculate
\allowdisplaybreaks{
\begin{eqnarray} \nonumber
\lefteqn{ \int_0^1 \frac{d}{dt} E_{k, \omega} (\phi) dt } \\ \nonumber
& = & \frac{k+1}{V} \int_0^1 \int_M (- \phi - t \dotp) ( \op - (1-t) \ddbp)^k \wedge \op^{n-k} dt \\ \nonumber
&& \mbox{} - \frac{n-k}{V} \int_0^1 \int_M \dotp ( (\op - (1-t) \ddbp)^{k+1} - \op^{k+1}) \wedge \op^{n-k-1} dt \\ \label{eqn1}
& = & \frac{1}{V} \sum_{i=1}^k (k+1) \binom{k}{i} (-1)^{i+1} (A(i) + B(i)) + \frac{1}{V} \sum_{i=1}^{k+1} (n-k) \binom{k+1}{i} (-1)^{i+1} C(i), \qquad \quad
\end{eqnarray}}
where
\begin{eqnarray*}
A(i) & = & \int_0^1 \int_M (1-t)^i \phi (\ddbp)^i \wedge \op^{n-i} dt \\
B(i) & = & \int_0^1 \int_M t(1-t)^i \dotp (\ddbp)^i \wedge \op^{n-i} dt \\
C(i) & = & \int_0^1\int_M (1-t)^i \dotp (\ddbp)^i \wedge \op^{n-i} dt.
\end{eqnarray*}
Now,
\allowdisplaybreaks{
\begin{eqnarray}
C(i+1)  \label{eqn2}
& = & \frac{i+1}{n-i} (A(i) + B(i) - C(i)) - \frac{1}{n-i} D(i),
\end{eqnarray}}
where
$$D(i) = \left.  \int_M \phi (\ddbp)^i \wedge \op^{n-i} \right|_{t=0}.$$
Using the relation (\ref{eqn2}) in (\ref{eqn1}), we obtain
\allowdisplaybreaks{
\begin{eqnarray} \nonumber
\lefteqn{ \int_0^1 \frac{d}{dt} E_{k, \omega} (\phi) dt } \\ \nonumber
& = & \frac{1}{V} \sum_{i=1}^{k+1}  \left[ (n-i+1) \binom{k+1}{i} - (k+1) \binom{k}{i}
 - (n-k) \binom{k+1}{i}  \right] (-1)^i C(i)\\ \nonumber
 && \mbox{} + \frac{n(k+1)}{V} C(1)  + \frac{1}{V} \sum_{i=1}^{k} \binom{k+1}{i+1} (-1)^{i+1} D(i) \\ \label{eqn3}
 & = & \frac{n(k+1)}{V} C(1)+ \frac{1}{V} \sum_{i=1}^{k} \binom{k+1}{i+1} (-1)^{i+1} D(i),
\end{eqnarray}}
using the fact that
\begin{eqnarray} \label{eqnzero}
(n-i+1) \binom{k+1}{i} - (k+1) \binom{k}{i}
 - (n-k) \binom{k+1}{i} & = & 0.
 \end{eqnarray}
Now observe that
\begin{equation} \label{eqn4}
C(1)= \frac{1}{n} \int_0^1 \int_M (1-t) \phi (\Delta_{\phi} \dotp) \op^n dt,
\end{equation}
and
\allowdisplaybreaks{
\begin{eqnarray} \nonumber
\lefteqn{\sum_{i=1}^{k} \binom{k+1}{i+1} (-1)^{i+1} D(i)} \\ \nonumber
& = & \left. \sum_{i=1}^k \sum_{j=0}^{i-1} \binom{k+1}{i+1} \binom{i-1}{j} (-1)^{i+j}  \int_M  \sqrt{-1} \partial \phi \wedge \dbar \phi \wedge \omega^j \wedge \omega_{\phi}^{n-j-1} \right|_{t=0} \\ \nonumber
& = & \left.  \sum_{j=0}^{k-1} \left[  \sum_{i=j+1}^k \binom{k+1}{i+1} \binom{i-1}{j} (-1)^{i+j}  \right] \int_M  \sqrt{-1} \partial \phi \wedge \dbar \phi \wedge \omega^j \wedge \omega_{\phi}^{n-j-1} \right|_{t=0}\\ \label{eqn5}
& = &  \left.  \sum_{j=0}^{k-1} (j-k) \int_M  \sqrt{-1} \partial \phi \wedge \dbar \phi \wedge \omega^j \wedge \omega_{\phi}^{n-j-1} \right|_{t=0},
\end{eqnarray}}
where we are using the identity
$$\sum_{i=j+1}^k \binom{k+1}{i+1} \binom{i-1}{j} (-1)^{i+j} = j-k.$$
Combining (\ref{eqn3}), (\ref{eqn4}) and (\ref{eqn5}) completes the proof of lemma.  Q.E.D.
\end{proof}

We have proved the first part of Theorem \ref{theorem1} if
$\Rc(\tilde{\omega}) >0$.
For the semi-definite case $\Rc(\tilde{\omega}) \ge 0$ we perturb
$\tilde{\omega}$ to $\tilde{\omega}_{\alpha} \in [\tilde{\omega}]$
for small $\alpha>0$ using the equation $\tilde{\omega}_{\alpha}^n =
e^{\alpha \tilde{f}} \tilde{\omega}^n$ for $\tilde{f}$ defined by
$\Rc({\tilde{\omega}}) - \tilde{\omega} = \ddbar \tilde{f}$ and
normalized appropriately.  Notice that
$\Rc(\tilde{\omega}_{\alpha})>0$.  We can then apply the argument
above to obtain $E_k(\omega_{KE}, \tilde{\omega}_{\alpha}) \ge 0$.
Let $\alpha$ tend to zero and observe that by Yau's estimates
\cite{Ya1}, $\tilde{\omega}_{\alpha} \rightarrow \tilde{\omega}$ in
$C^{\infty}$.    Hence $E_k (\omega_{KE}, \tilde{\omega}) \ge 0$.

For the second part of Theorem \ref{theorem1},  assume that $E_k
(\omega_{KE},\tilde{\omega}) = 0$ with $\Rc(\tilde{\omega}) \ge 0$.
We will show that $\tilde{\omega}$ is K\"ahler-Einstein.  We can
suppose $k \ge 1$ since the result is already known for $k=0$.  With
$\tilde{\omega}_{\alpha}$ as above, set $\omega_{\alpha} = \Rc
(\tilde{\omega}_{\alpha}) >0$.
 For each $\alpha>0$, by Lemma \ref{lemmaBM}, there exist families of metrics $\omega^{\epsilon}_{\alpha} $ and $\tilde{\omega}^{\epsilon}_{\alpha} = \omega^{\epsilon}_{\alpha} + \ddbar \psi^{\epsilon}_{\alpha, 0}$ with $\omega^{\epsilon}_{\alpha} \rightarrow \omega_{\alpha}$ and $\tilde{\omega}^{\epsilon}_{\alpha} \rightarrow \tilde{\omega}_{\alpha}$ in $C^{\infty}$ as $\epsilon \rightarrow 0$, and such that $\psi^{\epsilon}_{\alpha, t}$ satisfies
$$(\omega^{\epsilon}_{\alpha} + \ddbar \psi^{\epsilon}_{\alpha, t})^n = e^{-t \psi^{\epsilon}_{\alpha, t} + f^{\epsilon}_{\alpha}} (\omega^{\epsilon}_{\alpha})^n,
 \quad \omega^{\epsilon}_{\alpha} + \ddbar \psi^{\epsilon}_{\alpha, t}>0 ,$$
for $t \in [0,1]$, where $\Rc(\omega^{\epsilon}_{\alpha}) -
\omega^{\epsilon}_{\alpha} = \ddbar f^{\epsilon}_{\alpha}$. From
Lemma \ref{lemmaEkformula1} we see that
\begin{eqnarray} \label{eqnfinal}
E_k (\omega_{KE}, \tilde{\omega}^{\epsilon}_{\alpha}) \ge
\frac{k}{V}  \int_M \sqrt{-1} \partial \psi^{\epsilon}_{\alpha, 0}
\wedge \overline{\partial}  \psi^{\epsilon}_{\alpha, 0} \wedge
(\tilde{\omega}_{\alpha}^{\epsilon})^{n-1}.
\end{eqnarray}
First let $\epsilon$ tend to zero in (\ref{eqnfinal}).  The
functions $\psi_{\alpha, 0}^{\epsilon}$, normalized to have zero
integral with respect to $\omega_{KE}$,
 converge smoothly to $\psi_{\alpha,0}$ satisfying $\tilde{\omega}_{\alpha} = \omega_{\alpha} + \ddbar \psi_{\alpha,0}$, and $E_k (\omega_{KE}, \tilde{\omega}^{\epsilon}_{\alpha})$ tends to $E_k (\omega_{KE}, \tilde{\omega}_{\alpha})$.  Now let $\alpha$ tend to zero.  By Yau's estimates again, $\psi_{\alpha,0}$ converges in $C^{\infty}$ to $\psi_0$ satisfying $\tilde{\omega} = \Rc(\tilde{\omega}) + \ddbar \psi_0$.  Since $E_k (\omega_{KE}, \tilde{\omega})=0$, the inequality (\ref{eqnfinal}) becomes
\begin{eqnarray} \nonumber
0 \ge \frac{k}{V}  \int_M \sqrt{-1} \partial \psi_0\wedge
\overline{\partial}  \psi_0 \wedge \tilde{\omega}^{n-1}.
\end{eqnarray}
Hence $\psi_0$ is a constant and  $\tilde{\omega} =
\Rc(\tilde{\omega})$.  Q.E.D.

\bigskip
\bigskip
\pagebreak[3]
\setcounter{equation}{0}
\setcounter{lemma}{0}
\addtocounter{section}{1}
\noindent
{\bf 4. The lower bound of $E_1$}
\bigskip
\nopagebreak

In this section we give a proof of Theorem \ref{theorem2}.  Let us first deal with the easy case when $c_1(M)=0$.   From the formula obtained in section 2, we see that
$$E_{1,\omega_{KE}} (\phi) = \frac{1}{V} \int_M \log \left( \frac{\omega_{\phi}^n}{\omega^n} \right) \wedge  \Rc(\omega_{\phi}) \wedge \op^{n-1}.$$
Using the definition of $\Rc(\op)$ and integration by parts, we obtain:
\begin{eqnarray*}
E_{1, \omega_{KE}} (\phi) & = &\frac{1}{V} \int_M \sqrt{-1} \partial \log \left( \frac{\omega_{\phi}^n}{\omega^n} \right) \wedge \dbar \log \left(  \frac{\omega_{\phi}^n}{\omega^n} \right) \wedge \omega_{\phi}^{n-1} \ge  0,
\end{eqnarray*}
as required.

We turn now to the case when $c_1(M)>0$.  Consider the equation
\begin{equation} \label{eqnma2}
\omega_{\psi_t}^n = e^{t f + c_t} \omega^n,
\end{equation}
where $f$ satisfies
$\Rc(\omega) - \omega = \ddbar f$ and $\int_M e^f \omega^n = V,$
and $c_t$ is the constant chosen so that
$\int_M e^{t f + c_t} \omega^n = V.$
By the theorem of Yau \cite{Ya1}, we know that there is a unique $\psi_t$ with $\int_M \psi_t \omega^n =0$ solving (\ref{eqnma2}) for $t$ in $[0,1]$.  Notice that $\psi_0 = 0$.
Differentiating the equation we obtain
\begin{equation} \label{eqnddtpsi}
\Delta_{\psi_t} \dotpsi_t = f + \textrm{constant},
\end{equation}
and
\begin{equation} \label{eqnRc2}
\Rc(\omega_{\psi_t}) = \omega_{\psi_t} + (1-t) \ddbar f - \ddbpsi_t.
\end{equation}
The following lemma is the key result of this section.

\begin{lemma} \label{lemmaEk2} Let $\psi=\psi_t$ be a solution of (\ref{eqnma2}).  Then
\begin{eqnarray} \nonumber
E_{k, \omega} (\psi_1)
& = & - \frac{1}{V} \sum_{i=0}^{k-1} \frac{(n-k)(i+1)}{n+1} \left. \int_M \sqrt{-1} \partial \psi \wedge \dbar \psi \wedge \omega^i \wedge \opsi^{n-i-1} \right|_{t=1} \\ \nonumber
&& \mbox{} \left. - \frac{1}{V} \sum_{i=k}^n (k+1) \frac{n-i}{n+1} \int_M \sqrt{-1} \partial \psi \wedge \dbar \psi \wedge \omega^i \wedge \opsi^{n-i-1} \right|_{t=1} \\ \nonumber
&& \mbox{} - \frac{k+1}{V} \int_0^1 \int_M (1-t) (\Delta_{\psi} \dot{\psi})^2 \, \opsi^n dt \\ \nonumber
\ \label{eqnEkpsi}
&& \mbox{} + \frac{1}{V} \sum_{i=1}^k  \binom{k+1}{i+1}  \int_M f \left( \ddbar f \right)^i \wedge \omega^{n-i}.
\end{eqnarray}
In particular,
$E_{1, \omega} (\psi_1) \le 0.$
\end{lemma}

Given this lemma, we will prove Theorem \ref{theorem2}.  Let $\omega' \in \kahler$ be given.  Then by Yau's theorem there exists $\tilde{\omega} \in \kahlerp$ such that $\Rc (\tilde{\omega}) = \omega'$.  By the cocycle condition and Theorem \ref{theorem1},
\begin{eqnarray*}
E_1 (\omega_{KE}, \omega') & = & E_1 (\omega_{KE}, \tilde{\omega}) + E_1 (\tilde{\omega}, \omega') \\
& \ge & E_1 (\tilde{\omega}, \omega').
\end{eqnarray*}
We apply Lemma \ref{lemmaEk2} with $\omega = \omega'$.  Then we see that  $\omega_{\psi_1} = \tilde{\omega}$, and
$$E_{1, \omega'}(\psi_1) = E_1 (\omega', \tilde{\omega}) \le 0,$$
thus completing the first part of the proof of Theorem \ref{theorem2}.  The second part, that equality implies K\"ahler-Einstein, will follow by a similar argument to that given in section 3.
It remains to prove Lemma \ref{lemmaEk2}.

\pagebreak[3]
\bigskip
\noindent
{\bf Proof of Lemma \ref{lemmaEk2}}
\bigskip

Calculate using (\ref{eqnddtpsi}) and (\ref{eqnRc2}):
\begin{eqnarray} \nonumber
\lefteqn{
E_{k, \omega} (\psi_1)  } \\ \nonumber
& =  & \frac{k+1}{V} \int_0^1 \int_M (\Delta_{\psi} \dotpsi) (\opsi + (1-t) \ddbar f - \ddbpsi)^k \wedge \opsi^{n-k} dt \\ \nonumber
&& \mbox{} - \frac{n-k}{V} \int_0^1 \int_M \dotpsi \left( ( \opsi + (1-t) \ddbar f - \ddbpsi)^{k+1} - \opsi^{k+1} \right) \wedge \opsi^{n-k-1} dt \\ \label{eqn10}
& = & \frac{1}{V} \sum_{i=1}^k (k+1) \binom{k}{i} \hat{A}(i) - \frac{1}{V} \sum_{i=1}^{k+1} (n-k) \binom{k+1}{i} \hat{B}(i),
\end{eqnarray}
where
$$\hat{A}(i) = \int_0^1 \int_M f \left( (1-t) \ddbar f - \ddbpsi\right)^i \wedge \opsi^{n-i} dt,$$
and
$$\hat{B}(i) = \int_0^1 \int_M \dotpsi \left( (1-t) \ddbar f - \ddbpsi \right)^i \wedge \opsi^{n-i}dt.$$
Now calculate
\begin{eqnarray} \label{eqn9}
\hat{B}(i+1)
 =  \frac{i+1}{n-i} \hat{A}(i) + \frac{i+1}{n-i} \hat{B}(i) - \frac{1}{n-i} ((-1)^i\hat{C}(i) + \hat{D}(i)),
\end{eqnarray}
where
$$\left. \hat{C}(i) = \int_M \psi \left( \ddbpsi \right)^i \wedge \opsi^{n-i} \right|_{t=1} \ \  \textrm{and}
\ \  \hat{D}(i) = \int_M f \left( \ddbar f \right)^i \wedge \omega^{n-i}.$$
Making use of the relation (\ref{eqn9}) in (\ref{eqn10}), as well as (\ref{eqnzero}), we see that
\begin{eqnarray} \label{eqn11}
E_{k, \omega} (\psi_1)
& = & - \frac{n(k+1)}{V} \hat{B}(1) + \frac{1}{V} \sum_{i=1}^k \binom{k+1}{i+1} ( (-1)^i \hat{C}(i) + \hat{D}(i)).
\end{eqnarray}
Now
\begin{equation} \label{eqnB1}
\hat{B}(1) = \frac{1}{n} \int_0^1 \int_M (1-t) (\Delta_{\psi} \dot{\psi})^2\,  \opsi^n dt - \frac{1}{n} \int_0^1\int_M \dotpsi (\Delta_{\psi} \psi)  \opsi^n dt.
\end{equation}
Using (\ref{eqnddtIminusJ}) we have
\begin{eqnarray}
\frac{1}{n} \int_0^1\int_M \dotpsi (\Delta_{\psi} \psi)  \opsi^n dt =  -\frac{V}{n} \int_0^1 \frac{d}{dt} (I_{\omega} - J_{\omega})(\psi) dt
 =  - \frac{V}{n} (I_{\omega} - J_{\omega})(\psi_1). \label{eqn12}
\end{eqnarray}
Let us put this all together.  By the same calculation as in (\ref{eqn5}), we have
$$\left. \frac{1}{V} \sum_{i=1}^k \binom{k+1}{i+1} (-1)^i \hat{C}(i) = \frac{1}{V} \sum_{i=0}^{k-1} (k-i) \int_M \sqrt{-1} \partial \psi \wedge \dbar \psi \wedge \omega^{i} \wedge \omega_{\psi}^{n-i-1} \right|_{t=1}.$$
Combining this with equations (\ref{eqn11}), (\ref{eqnB1}) and (\ref{eqn12}),
\begin{eqnarray} \nonumber
E_{k,\omega}(\psi_1) & = &  -(k+1) (I_{\omega} - J_{\omega})(\psi_1) \\ \nonumber
&& \mbox{} - \frac{(k+1)}{V} \int_0^1 \int_M (1-t) (\Delta_{\psi} \dot{\psi})^2 \omega_{\psi}^n dt \\ \nonumber
&& \mbox{} \left. + \frac{1}{V} \sum_{i=0}^{k-1} (k-i) \int_M \sqrt{-1} \partial \psi \wedge \dbar \psi \wedge \omega^{i} \wedge \omega_{\psi}^{n-i-1} \right|_{t=1} \\ \label{eqnEkpsi1formula}
&& \mbox{} +  \frac{1}{V} \sum_{i=1}^k  \binom{k+1}{i+1}  \int_M f \left( \ddbar f \right)^i \wedge \omega^{n-i}.
\end{eqnarray}
Combining the first and third terms completes the first part of the lemma.
Finally, to show that $E_1(\psi_1) \le 0$, just observe that, by integration by parts, the last term is nonpositive for $k=1$.  Notice that this is the only step that fails for $k>1$. This completes the proof of the lemma.  Q.E.D.

\bigskip
\bigskip
\setcounter{equation}{0}
\setcounter{lemma}{0}
\addtocounter{section}{1}
\setcounter{theorem}{0}
\pagebreak[3]
\noindent
{\bf 5. Properness of $E_1$}
\bigskip

In this section we will prove Theorem \ref{theorem3} and Theorem \ref{theorem4}.  In fact we will assume that $M$ has no nontrivial holomorphic vector fields, which will prove only the first part of these theorems.  The proof of the second part in each case is almost identical and we refer the reader to \cite{Ti2}.

\bigskip
\noindent
{\bf Proof of Theorem \ref{theorem3}}
\bigskip

Let $\theta$ be in $P(M, \omega_{KE})$ and set
$\omega = \omega_{KE} + \ddbar \theta.$
We consider the family of Monge-Amp\`ere equations  (\ref{eqnma1}), $\omega_{\phi_t}^n = e^{-t \phi_t + f} \omega^n$, corresponding to this particular choice of $\omega$,
where $f$ is given by
$\Rc(\omega) - \omega = \ddbar f$ and  $\int_M e^f \omega^n = V.$
Since there are no nontrivial holomorphic vector fields, we can obtain a solution $\phi_t$ for $0 \le t  \le 1$ \cite{BaMb}.  Moreover, $\phi_1 = - \theta + \textrm{constant}$.

We have the following lemma.

\begin{lemma} \label{lemmaE1formula}
$$E_{1, \omega_{KE}} (\theta) \ge 2 \int_0^1 (I_{\omega} - J_{\omega}) (\phi_t) dt.$$
\end{lemma}

\begin{proof}
We will work in more generality and derive a formula for $E_{k, \omega_{KE}}(\theta)$.
Using the cocycle condition for $E_k$ and Lemma \ref{lemmaEkformula1} we have
\begin{eqnarray*}
E_{k, \omega_{KE}} (\theta) & = &  -E_{k, \omega} (\phi_0) - \frac{k+1}{V} \int_0^1 \int_M (1-t) \phi_t (\Delta_{\phi_t} \dot{\phi}_t)\omega_{\phi_t}^n dt \\
&& \mbox{} + \frac{1}{V} \sum_{i=0}^{k-1} (k-i) \int_M \sqrt{-1} \partial \phi_0 \wedge \dbar \phi_0 \wedge \omega^i \wedge \omega_{\phi_0}^{n-i-1},
\end{eqnarray*}
since $E_{k, \omega_{KE}} (\theta) = - E_{k, \omega} (\phi_1)$.  Now observe that, using (\ref{eqnddtIminusJ}), we have
\begin{eqnarray*}
\lefteqn{\frac{k+1}{V} \int_0^1 \int_M (1-t) \phi_t (\Delta_{\phi_t} \dot{\phi}_t ) \omega_{\phi_t}^n dt} \\
&= & \mbox{} - (k+1)  \int_0^1 (I_{\omega} - J_{\omega})(\phi_t) dt + (k+1) (I_{\omega} - J_{\omega})(\phi_0).
\end{eqnarray*}
Hence
\begin{eqnarray} \nonumber
E_{k, \omega_{KE}} (\theta) & = & - E_{k, \omega} (\phi_0) + (k+1)  \int_0^1 (I_{\omega} - J_{\omega})(\phi_t) dt - (k+1) (I_{\omega} - J_{\omega})(\phi_0) \\
&& \mbox{} + \frac{1}{V} \sum_{i=0}^{k-1} (k-i) \int_M \sqrt{-1} \partial \phi_0 \wedge \dbar \phi_0 \wedge \omega^i \wedge \omega_{\phi_0}^{n-i-1}. \label{eqn20}
\end{eqnarray}
Now, as in section 4, let $\psi_t$ be a solution of (\ref{eqnma2}), for $\omega$ as above:
$ \omega_{\psi_t}^n = e^{tf + c_t} \omega^n,$
with $c_t$ the appropriate constant.  Then observe that
$\psi_1 = \phi_0 + \textrm{constant}.$
Then making use of the calculation (\ref{eqnEkpsi1formula}) from section 4, we have
\begin{eqnarray} \nonumber
E_{k, \omega} (\phi_0) & = & -(k+1) (I_{\omega} - J_{\omega})(\phi_0) - \frac{(k+1)}{V} \int_0^1 \int_M (1-t) (\Delta_{\psi_t} \dot{\psi}_t)^2 \omega_{\psi_t}^n dt \\ \nonumber
&& \mbox{} + \frac{1}{V} \sum_{i=0}^{k-1} (k-i) \int_M \sqrt{-1} \partial \phi_0 \wedge \dbar \phi_0 \wedge \omega^{i} \wedge \omega_{\phi_0}^{n-i-1}  \\ \label{eqn21}
&& \mbox{} +  \frac{1}{V} \sum_{i=1}^k  \binom{k+1}{i+1}  \int_M f \left( \ddbar f \right)^i \wedge \omega^{n-i},
\end{eqnarray}
since $E_{k, \omega}(\phi_0) = E_{k, \omega} (\psi_1)$.
Combining (\ref{eqn20}) and (\ref{eqn21}) we obtain
\begin{eqnarray*}
E_{k, \omega_{KE}}(\theta) & = & (k+1) \int_0^1 (I_{\omega} - J_{\omega})(\phi_t) dt + \frac{(k+1)}{V} \int_0^1 \int_M (1-t) (\Delta_{\psi_t} \dot{\psi}_t)^2 \omega_{\psi_t}^n dt \\
&& \mbox{} - \frac{1}{V} \sum_{i=1}^k  \binom{k+1}{i+1}  \int_M f \left( \ddbar f \right)^i \wedge \omega^{n-i}.
\end{eqnarray*}
We see that for $k=1$, the last two terms are nonnegative, and this proves the lemma.  Q.E.D.
\end{proof}

Now observe that we can now apply the argument of \cite{Ti2}.  Recall that Tian proves the properness of the functional $F_{\omega}$ using the formula
$F_{\omega}(\theta) =  \int_0^1 (I_{\omega} - J_{\omega}) (\phi_t) dt,$
with $\phi_t$ as above.  In light of Lemma \ref{lemmaE1formula} we can directly apply his argument to obtain:

\begin{theorem} \label{theoremtian}
There exists $\delta=\delta(n)$ such that for every constant $K>0$ there exist positive constants $C_1$ and $C_2$ depending on $K$ such that
for all $\theta$ in $P(M, \omega_{KE})$ satifsying
\begin{equation} \label{eqncondition}
\emph{osc}_M(\theta) \le K (1+ J_{\omega_{KE}} (\theta)),
\end{equation}
we have
$$E_{1, \omega_{KE}} (\theta) \ge C_1 J_{\omega_{KE}}(\theta)^{\delta} - C_2.$$
Here, the oscillation $\emph{osc}_M$ is defined by
$\emph{osc}_M (\theta) = \sup_M (\theta) - \inf_M (\theta).$
\end{theorem}

We will now apply the argument of Tian and Zhu \cite{TiZh} to show that $E_{1, \omega_{KE}}$ is proper on the full space of potentials.  The argument for $E_1$ differs in only one step, but we will outline the whole argument here.

We need some lemmas.

\begin{lemma} \label{lemmaoscphitminusphi1}
For $\phi_t$ as above, there exists a constant $C_3$ depending only on $\omega_{KE}$ such that for $t \ge \frac{1}{2}$,
$$\emph{osc}_M (\phi_t - \phi_1) \le C_3 (1 + J_{\omega_{KE}} (\phi_t - \phi_1)).$$
\end{lemma}
\begin{proof}  We omit the proof, since it can be found in \cite{TiZh}. Q.E.D.
\end{proof}

We need the following lemma, which we state in more generality than is actually needed here.

\begin{lemma} \label{lemmat1t2}
Let $\phi_t$ be a solution of
$$\omega_{\phi_t}^n = e^{-t \phi_t + f} \omega^n.$$
Then for $0 \le t_1 \le t_2 \le 1$, we have
\begin{eqnarray*}
\lefteqn{E_{1, \omega} (\phi_{t_2}) - E_{1, \omega}(\phi_{t_1}) } \\
& = & -2(1-t_2) (I_{\omega} - J_{\omega})(\phi_{t_2}) + 2(1-t_1) (I_{\omega} - J_{\omega})(\phi_{t_1}) \\
&& \mbox{} - 2\int_{t_1}^{t_2} (I_{\omega} - J_{\omega})(\phi_t) dt
 \left. + \left( \frac{(1-t)^2}{V} \int_M \sqrt{-1} \partial \phi_{t} \wedge \dbar \phi_{t} \wedge \omega_{\phi_{t}}^{n-1} \right) \right|_{t=t_1}^{t=t_2}.
\end{eqnarray*}
\end{lemma}
\begin{proof}
Calculate, as in section 3, dropping the subscript $t$,
\allowdisplaybreaks{
\begin{eqnarray*}
\lefteqn{E_{1, \omega}(\phi_{t_2}) - E_{1, \omega}(\phi_{t_1}) } \\
& = & - \frac{2}{V} \int_{t_1}^{t_2} \int_M (1-t) (- \phi - t \dotp) \ddbp \wedge \op^{n-1} dt \\
&& \mbox{} - \frac{n-1}{V} \int_{t_1}^{t_2} \int_M \dotp ( (\op - (1-t) \ddbp)^{2} - \op^{2}) \wedge \op^{n-2} dt \\
& = & \frac{2}{V} \int_{t_1}^{t_2} \int_M (1-t) \phi \ddbp \wedge \op^{n-1} dt - \frac{2}{V} \int_{t_1}^{t_2} \int_M (1-t)^2 \dot{\phi} \ddbp \wedge \op^{n-1}dt \\
&& \mbox{}  +\frac{2n}{V} \int_{t_1}^{t_2} \int_M (1-t) \dot{\phi} \ddbp \wedge \op^{n-1} dt \\
&& \mbox{}  - \frac{1}{V} \int_{t_1}^{t_2} \int_M (1-t)^2 \phi \ddbp \wedge \frac{d}{dt} (\op^{n-1}) dt.
\end{eqnarray*}}
Integrating by parts in $t$, and making use of (\ref{eqnddtIminusJ}), we obtain
\begin{eqnarray*}
E_{1, \omega}(\phi_{t_2}) - E_{1, \omega}(\phi_{t_1})
& = & - 2 \int_{t_1}^{t_2} (1-t) \frac{d}{dt} (I_{\omega} - J_{\omega})(\phi) dt \\
&& \mbox{} \left. + \left( \frac{(1-t)^2}{V} \int_M \sqrt{-1} \partial \phi \wedge \dbar \phi \wedge \omega_{\phi}^{n-1} \right) \right|_{t=t_1}^{t=t_2},\end{eqnarray*}
and the lemma follows after integrating by parts in the first term. Q.E.D.
\end{proof}

We will use this lemma in the proof of the following:

\begin{lemma} \label{lemmaE1J}
$$E_{1, \omega_{KE}} (\phi_t - \phi_1) \le 2n (1-t) J_{\omega_{KE}} (\theta).$$
\end{lemma}
\begin{proof}
We use the cocycle condition for $E_1$ together with Lemma \ref{lemmat1t2} for $t_1=t$ and $t_2=1$ to obtain
\begin{eqnarray*}
E_{1,\omega_{KE}} (\phi_t - \phi_1)
 & = & E_{1, \omega} (\phi_t) - E_{1,\omega} (\phi_1) \\
& = & - 2(1-t) (I_{\omega} - J_{\omega})(\phi_t) + 2 \int_t^1 (I_{\omega} - J_{\omega})(\phi_s) ds \\
&& \mbox{} + \frac{(1-t)^2}{V} \int_M \sqrt{-1} \partial \phi_t \wedge \dbar \phi_t \wedge \omega_{\phi_t}^{n-1} \\
& \le & 2 \int_t^1 (I_{\omega} - J_{\omega})(\phi_s) ds \\
& \le & 2(1-t) (I_{\omega} - J_{\omega})(\phi_1) \\
& \le & 2n (1-t)J_{\omega_{KE}}(\theta),
\end{eqnarray*}
where we have used the fact (\ref{eqnintegralinequality}) that $(I_{\omega} -J_{\omega}) (\phi_t)$ is increasing in $t$.  Q.E.D. \end{proof}

Now note from Lemma \ref{lemmaoscphitminusphi1} that $\phi_t - \phi_1$ satisfies the condition (\ref{eqncondition}) for $K=C_3$ and so we can apply Theorem \ref{theoremtian} to obtain
$$E_{1,\omega_{KE}} (\phi_t - \phi_1) \ge C_4 (J_{\omega_{KE}} (\phi_t - \phi_1))^{\delta} -C_5,$$
for $t \ge \frac{1}{2}$ and for some constants $C_4$ and $C_5$ depending only on $\omega_{KE}$.  Applying Lemma \ref{lemmaoscphitminusphi1} and Lemma \ref{lemmaE1J}, we see that for $t \ge \frac{1}{2}$,
$$2n (1-t) J_{\omega_{KE}} (\theta) \ge C_6 (\textrm{osc}_M (\phi_t - \phi_1))^{\delta} - C_7,$$
for constants $C_6$ and $C_7$ depending only on $\omega_{KE}$.

We need another lemma.

\begin{lemma} For any $t \in [0,1]$,
$$\int_0^1 (I_{\omega} - J_{\omega})(\phi_s) ds \ge (1-t) (I_{\omega} - J_{\omega}) (\phi_1) - 2n(1-t) \emph{osc}_M (\phi_t - \phi_1).$$
\end{lemma}
\begin{proof}  The proof can be found in \cite{Ti2}.  Q.E.D.
\end{proof}

We can now complete the proof of Theorem \ref{theorem3}.   We have for $t \ge \frac{1}{2}$,
\begin{eqnarray*}
E_{1, \omega_{KE}} (\theta) & \ge &  2 \int_0^1 (I_{\omega} - J_{\omega}) (\phi_t) dt \\
& \ge & 2(1-t) (I_{\omega} - J_{\omega}) (\phi_1) - 4 n(1-t) \textrm{osc}_M (\phi_t - \phi_1) \\
& \ge & \frac{2(1-t)}{n} J_{\omega_{KE}} (\theta) - 4n(1-t)C_6^{-1/\delta} (2n(1-t) J_{\omega_{KE}}(\theta) + C_7)^{1/\delta}.
\end{eqnarray*}
Setting
$(1-t) = \frac{1}{C_8} (1+ J_{\omega_{KE}}(\theta))^{-1+\delta},$
for $C_8 >>  C^{-1}_6$ completes the proof.  Q.E.D.

\bigskip
\bigskip \pagebreak[3]
\noindent
{\bf Proof of Theorem \ref{theorem4}}
\bigskip

We have to show that if $E_1$ is proper then there exists a K\"ahler-Einstein metric on $M$.  Fix a metric $\omega$ and consider again the family of Monge-Amp\`ere equations
$\omega_{\phi_t}^n = e^{-t \phi_t + f} \omega^n,$
where $f$ is given by
$\Rc(\omega) - \omega = \ddbar f$ and $\int_M e^f \omega^n = V.$
If we have a solution for $t=1$ then $\omega + \ddbar \phi_1$ is a K\"ahler-Einstein metric.
We use the continuity method.   We have a solution $\phi_0$ at $t=0$ by Yau's theorem \cite{Ya1}.
 It is well-known that we can find a solution $\phi_t$ for $0 \le t \le 1$ as long as we can bound the $C^0$ norm of $\phi_t$ uniformly in $t$.  Moreover,
$\| \phi_t \|_{C^0} \le \tilde{C} (1+ J_{\omega} (\phi_t))$
for a constant $\tilde{C}$ depending only on $\omega$.  So it suffices to show that $J_{\omega}(\phi_t)$ is uniformly bounded from above.  Since $E_{1, \omega}$ is proper, it suffices to show that $E_{1,\omega} (\phi_t)$ is bounded from above uniformly in $t$.

We apply Lemma \ref{lemmat1t2} in the case when $t_2=t$ and $t_1=0$.  Then
\begin{eqnarray*}
E_{1, \omega} (\phi_t) - E_{1, \omega}(\phi_0) & = & -2(1-t) (I_{\omega}- J_{\omega})(\phi_t) + 2(I_{\omega} - J_{\omega})(\phi_0) \\
&& \mbox{}
- 2 \int_0^t (I_{\omega} - J_{\omega})(\phi_s) ds \\ && \mbox{}
 + \frac{(1-t)^2}{V} \int_M \sqrt{-1} \partial \phi_t \wedge \dbar \phi_t \wedge \omega_{\phi_t}^{n-1} \\
&& \mbox{} - \frac{1}{V} \int_M \sqrt{-1} \partial \phi_0 \wedge \dbar \phi_0 \wedge \omega_{\phi_0}^{n-1} \\
& \le & 2(I_{\omega} - J_{\omega})(\phi_0) - 2 \int_0^t (I_{\omega} - J_{\omega})(\phi_s) ds \\ && \mbox{}
  - \frac{1}{V} \int_M \sqrt{-1} \partial \phi_0 \wedge \dbar \phi_0 \wedge \omega_{\phi_0}^{n-1},
\end{eqnarray*}
where we are using the definition of $I_{\omega} - J_{\omega}$.  It then follows immediately that $E_{1, \omega}(\phi_t)$ is bounded from above by a constant independent of $t$. Q.E.D.

\bigskip
\noindent {\bf Acknowledgements.} The authors would like to thank:
Professor D.H. Phong, for his constant support and advice; Professor
J. Sturm, for a number of very enlightening discussions; Professor
X.X. Chen, for some helpful suggestions and in particular for
bringing our attention to the functional $E_1$; Professor S.-T. Yau,
for his encouragement and support; Professor S.T. Paul for some
useful discussions; and the referees for some helpful comments.

\let\oldbibliography\thebibliography
\renewcommand{\thebibliography}[1]{%
  \oldbibliography{#1}%
  \setlength{\itemsep}{0pt}%
}


\begin{thebibliography}{99} \small
\bibitem[Au]{Au} Aubin, T. {\em \'Equations du type Monge-Amp\`ere sur les
vari\'et\'es k\"ahleriennes compactes}, Bull. Sci. Math. (2) {\bf 102} (1978),
no. 1, 63--95
\bibitem[Ba]{Ba} Bando, S. {\em The K-energy map, almost Einstein K\"ahler
metrics and an inequality of the Miyaoka-Yau type}, Tohuku Math. Journ. {\bf 39}
(1987), 231--235
\bibitem[BaMb]{BaMb} Bando, S. and Mabuchi, T. {\em Uniqueness of Einstein
K\"ahler metrics modulo connected group actions}, Adv. Stud. in Pure Math. {\bf 10}
(1987), 11--40
\bibitem[Ch1]{Ch1} Chen, X.X. {\em On the lower bound of the Mabuchi energy and
its application}, Int. Math. Res. Not. {\bf 12} (2000), 607--623
\bibitem[Ch2]{Ch2} Chen, X.X. {\em On the lower bound of energy functional $E_1$ (I)-- a stability theorem on  the Kaehler Ricci flow}, preprint, math.DG/050219
\bibitem[Ch3]{Ch3} Chen, X.X., {\em On the lower bound of energy functional $E_1$ (II)}, unpublished
\bibitem[ChTi1]{ChTi1} Chen, X.X. and Tian, G. {\em Ricci flow on
K\"ahler-Einstein surfaces}, Invent. Math. {\bf 147} (2002), 487--544
\bibitem[ChTi2]{ChTi2} Chen, X.X. and Tian, G. {\em Ricci flow on
K\"ahler-Einstein manifolds}, Duke Math. J. {\bf 131} (2006), no. 1, 17--73
\bibitem[ChTi3]{ChTi3} Chen, X.X. and Tian, G. {\em Geometry of K\"ahler metrics and foliations by holomorphic discs}, preprint, math.DG/0507148
\bibitem[De]{De} Deligne, P. {\em Le determinant de la cohomologie}, Contemporary Math. {\bf 67} (1987), 93-177
\bibitem[Di]{Di} Ding, W. {\em
Remarks on the existence problem of positive K\"ahler-Einstein metrics},
Math. Ann. {\bf 282} (1988), no. 3, 463--471.
\bibitem[DiTi]{DiTi} Ding, W. and Tian, G. {\em K\"ahler-Einstein metrics and
the generalized Futaki invariant}, Invent. Math. {\bf 110} (1992), 315-335
\bibitem[Do1]{Do1} Donaldson, S. K. {\em Scalar curvature and projective
embeddings, I.}, J. Diff. Geom. {\bf 59} (2001), no. 3, 479-522
\bibitem[Do2]{Do2} Donaldson, S.K. {\em Scalar curvature and stability of toric
varieties}, J. Diff. Geom. {\bf 62} (2002), no. 2, 289--349
\bibitem[Do3]{Do3} Donaldson, S.K. {\em Scalar curvature and projective embeddings, II}, Q. J. Math. {\bf 56} (2005), no. 3, 345--356
\bibitem[Fu]{Fu} Futaki, A. {\em An obstruction to the existence of Einstein
K\"ahler metrics}, Invent. Math. {\bf 73} (1983), 437--443
\bibitem[HwMs]{HwMs} Hwang, A.D. and Maschler, G. {\em Central K\"ahler metrics
with non-constant central curvature}, Trans. Amer. Math. Soc. {\bf 355} (2003),
no. 6, 2183--2203
\bibitem[Mb]{Mb} Mabuchi, T. {\em K-energy maps integrating Futaki invariants},
T\^{o}hoku Math. Journ. {\bf 38} (1986), 575--593
\bibitem[Ms]{Ms} Maschler, G. {\em Central K\"ahler metrics}, Trans. Amer. Math. Soc. {\bf 355} (2003), 2161--2182
\bibitem[PaTi]{PaTi} Paul, S. and Tian, G. {\em Algebraic and analytic K-stability},
Int. Math. Res. Not. {\bf 48} (2004), 2555--2591
\bibitem[PhSt1]{PhSt1} Phong, D.H. and Sturm, J. {\em Stability, energy
functionals, and K\"ahler-Einstein metrics}, Comm. Anal.
Geom. {\bf 11} (2003), no. 3, 565--597
\bibitem[PhSt2]{PhSt2} Phong, D.H. and Sturm, J. {\em Scalar curvature, moment
maps and the Deligne Pairing}, Amer. J. Math. {\bf 126} (2004), no. 1-2, 693-712
\bibitem[PhSt3]{PhSt3} Phong, D.H. and Sturm, J. {\em The Futaki invariant and the Mabuchi energy of a complete intersection}, Comm. Anal. Geom. {\bf 12} (2004), no. 1-2, 323--343
\bibitem[PhSt4]{PhSt4} Phong, D.H. and Sturm, J. {\em On stability and the convergence of the K\"ahler-Ricci flow}, J. Diff. Geom. {\bf 72} (2006),
no. 1, 149--168
\bibitem[RoTh]{RoTh} Ross, J. and Thomas, R. {\em A study of the Hilbert-Mumford criterion for the stability of projective varieties}, preprint, math.AG/0412519
\bibitem[Si]{Si} Siu, Y-T. {\em Lectures on Hermitian-Einstein metrics for
stable  bundles and K\"{a}hler-Einstein metrics}, Birkh\"{a}user Verlag,
Basel 1987
\bibitem[SoWe]{SoWe} Song, J. and Weinkove, B. {\em On the convergence and singularities of the J-flow with applications to the Mabuchi energy}, preprint, math.DG/0410418
\bibitem[Ti1]{Ti1} Tian, G. {\em The K-energy on hypersurfaces and stability},
Comm. Anal. Geom. {\bf 2} (1994), no. 2, 239--265
\bibitem[Ti2]{Ti2} Tian, G. {\em K\"ahler-Einstein metrics with positive scalar
curvature}, Invent. Math. {\bf 137} (1997), 1--37
\bibitem[Ti3]{Ti3} Tian, G. {\em Canonical metrics in K\"ahler geometry},
Lectures in Mathematics, ETH Z\"urich, Birkhauser Verlag, Basel 2000
\bibitem[TiZh]{TiZh} Tian, G. and Zhu, X. {\em A nonlinear inequality of Moser-Trudinger type}, Calc. Var. {\bf 10} (2000), 349--354
\bibitem[We]{We} Weinkove, B. {\em
On the $J$-flow in higher
dimensions and the lower boundedness of the
Mabuchi energy},  J. Diff. Geom. {\bf 73} (2006), no. 2, 351--358
\bibitem[Ya1]{Ya1} Yau, S.-T. {\em On the Ricci curvature of a compact K\"ahler
manifold and the complex Monge-Amp\`ere equation, I}, Comm. Pure Appl. Math. {\bf 31}
(1978), 339--411
\bibitem[Ya2]{Ya2} Yau, S.-T. {\em Open problems in geometry}, Proc. Symposia Pure
Math. {\bf 54} (1993), 1--28 (problem 65)
\bibitem[Zh]{Zh} Zhang, S. {\em Heights and reductions of semi-stable varieties}, Compositio Math. {\bf 104} (1996), no. 1, 77--105
\end{thebibliography}
\end{document}